\documentclass[Proceedings,letterpaper]{ascelike-new}
\usepackage[utf8]{inputenc}
\usepackage[T1]{fontenc}
\usepackage{lmodern}
\usepackage{graphicx}
\usepackage[figurename=Fig.,labelfont=bf,labelsep=period]{caption}
\usepackage{subcaption}
\usepackage{amsmath}
\usepackage{newtxtext,newtxmath}
\usepackage[ruled,vlined]{algorithm2e}
\usepackage[colorlinks=false,citecolor=red,linkcolor=black]{hyperref}
\NameTag{Patil, \today}
\begin{document}

% You will need to make the title all-caps
\title{Budget-constrained rail electrification modeling using symmetric traffic assignment - a North American case study}

\author[1]{Priyadarshan Patil}
\author[2]{Rydell Walthall}
\author[2]{Stephen D. Boyles, Ph.D.}

\affil[1]{Graduate Program in Operations Research and Industrial Engineering, The University of Texas at Austin, USA. Email: Priyadarshan@utexas.edu}
\affil[2]{Department of Civil, Architectural and Environmental Engineering, The University of Texas at Austin, USA.}

\maketitle

\begin{abstract}
We consider a budget constrained rail network electrification problem with associated changes in costs of energy usage (via path gradient and curvature), operations, and long-term maintenance. In particular, we consider that freight flows on such a network form a user equilibrium.  Interactions between electric and diesel trains on the same corridor are represented with nonseparable link performance functions, which nevertheless have a symmetric Jacobian. This bi-level formulation is solved for the North American railroad network using a genetic algorithm (GA), incorporating domain-specific insights to reduce the number of solutions which must be considered.  We analyze solution characteristics and decision-making implications. Results show that broad connectivity would be beneficial for most impact. Increasing demand shifts electrified corridors towards the more populous east and gulf coasts, while increased operational costs results in electrification of routes through mountainous terrains.
\end{abstract}

\section{Introduction}
\label{sec:intro}
Rail networks play a vital role in local and national economic structures. In many countries such as Japan and Switzerland, they constitute a significant share of passenger transport mode share. Rail freight transit also accounts for a large portion of total freight transit, exceeding 50\% modal share in large economies like Canada and Russia. Therefore, a significant opportunity exists for government policies incentivizing rail network improvements to provide for larger economic, environmental and social returns.

Rail network electrification, and the accompanying transition from diesel-electric to fully electric locomotives, are important steps towards sustainable systems and renewable fuel sources for the United States of America.  There are many studies on the impact and cost-benefit analysis of rail electrification~\cite{USDOT,USDOT2}. Advantages of electric locomotion include lower long-term energy and locomotive maintenance costs, lower noise and air pollution levels, faster acceleration, and more flexibility in the primary power source, leading to less volatility from fuel price fluctuations. These benefits must be balanced with significant upfront investment for infrastructure upgrades, higher infrastructure maintenance costs, vulnerability of overhead architecture, higher property tax obligations for private rail companies, and the general uncertainty in the investment~\cite{walthall}.

We consider a rail network design problem (RNDP) where parts of a freight rail network can be electrified, subject primarily to budget constraints. The RNDP is formulated as a bi-level problem: the upper level problem deciding optimal subset of links for electrification, and the lower level problem calculating link flows as well as associated network metrics. Our upper-level formulation uses the objective of minimizing private costs, laying a foundation that can be adapted to improve net social benefit and reflect where subsidies ought to be directed.  Given that there are many rail operators, and shippers can choose which operator to use (not necessarily aligning with net social benefit), we model the lower level problem as a user-equilibrium traffic assignment problem (TAP).  Genetic algorithms (GA) solve the higher level problem of finding the optimal links to electrify.

In this setting, rather than apportioning an electrification budget to each rail operator independently, a utilitarian schema allocates the electrification budget for specific link electrification in order to bring about the greatest possible cost reductions across the network. The rail operators and shippers then respond to these changes by altering their scheduling and flow patterns to minimize their individual costs. Given multiple operators, the lower level problem is a setting where flow is directed ``selfishly,'' to minimize shipment costs. For the lower level problem, we assume that these activities lead to an equilibrium, where the costs of shipment flows cannot be lowered unilaterally. With shipment flow expressed in tons, as a continuous quantity over the long term, this problem satisfies the traffic assignment user equilibrium assumptions. Given few self-optimizing fleet owners, we have a Nash-Cournot equilibrium, which in the limit results in user equilibrium flow pattern. \citeN{van1991multiple} show that a Nash-Cournot equilibrium for two fleets results in less than 5\% difference in avg. travel times (and by extension, total system travel time) for large networks. As the number of fleets increase, this difference decreases and approaches zero. Therefore, we use the user equilibrium assumption in our study. This assumption is consistent with prior literature on the topic ~\cite{uddin2015,wang2018modeling}.

\subsection{Contributions and overview}
\label{sec:contributions_overview}

NDPs are widely studied in road networks.  Rail NDPs vary in two significant ways.  First, most of the freight rail network is privately owned by the user or contracted out for usage, which leads to non-socially-optimal usage restrictions.
Rail electrification has highly uncertain, and possibly negative, rates of return when external benefits are not accounted for.  This paper provides a framework for future analysis of policy interventions in the US to internalize the benefits to the private companies that own the rails and would be responsible for implementing electrification. Second, the characteristics of individual network links (track curvature and gradient) affect maintenance and operating costs more so than in road networks.

\noindent With these distinctions in mind, the main contributions of this article are as follows:
\begin{itemize}
    \item We formulate the rail electrification NDP as a bi-level optimization problem incorporating electrification costs; fuel, locomotive, and operational costs; and train resistance (bearing, flange, air, grade, curve, braking, and inertia) costs.
    \item We derive the appropriate flow shift formula for Algorithm B for solving our formulation (traffic assignment with symmetric link interactions) and show that it meets the optimality conditions required for convergence
    \item We use a general-purpose metaheuristic (a genetic algorithm) based on problem-specific insights to generate high-quality solutions, and solve this problem on a large-scale North American rail network.
    \item We conduct sensitivity analysis and analyze the resulting solutions to draw insights and policy conclusions.
\end{itemize}

The rest of this article is organized as follows. We first provide background information on rail electrification and discuss advances in traffic assignment, NDPs, and solution methods for both problems. We then describe the formulation for the rail electrification NDP and associated model components. We then describe the North American rail network dataset and demand data we use, and outline our experiment design.  We follow this with a summary of the results from our experiments and draw practical insights.  We conclude by summarizing our findings and suggest avenues for future work.

%_________________________________________________________________________

\section{Literature review}
\label{sec:lit_review}

\subsection{Rail electrification in the US}

The oil crisis of the 1970's spurred research into the economics of rail electrification in the US.  In 1977, \citeN{schwarm1977} examined the detailed costs associated with rail electrification infrastructure, finding three primary categories: power delivery, public works compatibility, and signaling systems compatibility. Power delivery uas used to calculate electrified link capacity. \citeN{Kneschke1986} then detailed the design requirements of the traction substations, which affect the electrified lines' capacities.  Meanwhile, other researchers looked at benefits of electrification. \citeN{Whitford1981} examined the energy savings from electrifying the high-density portions of the rail network, and \citeN{Ditmeyer1985} incorporated ancillary benefits and costs, such as those accruing to maintenance, reliability, and fuel handling.  At the time, emissions were not a major concern.  Additionally, computational costs and the relative youth of network optimization techniques meant that these analyses primarily used traffic density as a heuristic for link selection.

As computational resources and optimization algorithms progressed, the oil crisis abated and deregulation led to a decrease in rail network mileage and interest in rail electrification.  In 2012, \citeN{Cambridgesystematics2012} conducted an extensive study on electrification for the Southern California Association of Governments. \citeN{RailTEC2016} conducted a similar study for the California Air Resource Board in 2016. Both studies were primarily motivated by reducing rail emissions in dense urban areas. The latter study highlighted the problem of network connectivity: because rail freight routes are particularly long, and the costs of switching from electric to diesel-electric locomotives en route are high (in terms of delay and logistical costs of ensuring locomotives are available), and confining electrification to relatively small regions is not an efficient way to limit emissions. Outside of regulatory constraints or internalized emissions costs, rail companies would likely send large numbers of diesel-electric trains through electrified portions of the rail network in order to limit overall route costs. We thus face an optimization problem impacted by policy: how should the initial links be selected for electrification, given budget constraints and a desire to reduce emissions from the rail network?

There are very few studies to draw on for link improvement under a budget constraint, such as \citeN{Mishra2016}, which applied numerical methods towards ``[making] optimal investment decisions... in moderate and large transportation networks.''  Mishra utilized travel costs to road users, primarily consisting of temporal costs.  Applying a similar framework to a rail network considers similar costs, although fuel costs for rail constitute a relatively larger portion, and time costs a relatively smaller one.

\subsection{Rail electrification of other national networks}

While the European rail network is extensively electrified,  most of it is used for passenger transport and comparisons to the US network are difficult. This can be attributed to deliberate prioritization of passenger rail efficiency over freight, as evidenced by limits on maximum permissibly train lengths, maximum axle loads, and vertical car height \cite{walthall}. Similarly, the Japanese rail network is used primarily for passengers, accounting for less than 1\% of total freight load.

India and China both offer examples of extensive electrification of freight rail systems relevant to the US network. China has extensively electrified both its passenger and freight networks in the past decades, with over 70\% total electrification and about 25\% of its network electrified between 1990 and 2007~\cite{chineserail}. Much of this electrification has coincided with major construction to upgrade capacity and signaling infrastructure, especially in mountainous areas.

Similarly, India has extensively electrified its network (about 70\%), prioritizing selected main lines and high density routes~\cite{indianrail}. India has constructed high overhead contact systems in order to accommodate double-stacked freight.  Similar construction would be required in the US, where double-stacked container freight is a large component of American railroads’ revenues.

\subsection{Traffic assignment problem background}

This lower level assignment problem is a generalized version of the static traffic assignment problem formulated as a convex program by \citeN{Beckmann}. Specifically, the rail traffic assignment can be modeled as a variant of static TAP with symmetric link interactions if the flows refer to the tonnage of goods~\cite{SAM,uddin2015,wang2018modeling}.  We therefore assume the units of flow to be tons, allowing us to use existing solution methods to solve this problem. Given this formulation, the assignment problem follows Wardrop's first principle\cite{Wardrop}: ``All used paths between each origin and destination must have equal and minimal cost.'' Static TAP is a well-studied problem; for more on its properties and solution methods, see the books by \citeN{Patriksson} and \citeN{blubook_vol1_v085}. Static TAP solution methods include specialized algorithms to obtain the optimal solution. Some widely used methods are gradient projection~\cite{jaya}, Algorithm B~\cite{Dial}, and TAPAS~\cite{bar2010}, among others. A recent computational comparison by \citeN{perederieieva2015framework} concludes that Algorithm B and TAPAS outperform other tested methods. 

However, solution existence and optimality for these methods has been shown only for TAP with separable link costs, where the cost of each link depends only on its own flow. Generalizations have been explored for cases with link interactions, i.e., the cost of a link depends on the flow of many links.  The convex programming formulation is valid for symmetric interactions, i.e., when the Jacobian of the cost mapping is symmetric everywhere.  This extension to TAP was first mentioned by \citeN{prager1954problems}, discussing the need for traffic interactions on a two-way street. \citeN{dafermos1971extended} and \citeN{dafermos1972traffic} first modeled symmetric interactions using the entire flow vector and showed multi-class TAP and symmetric TAP as equivalent. Dafermos also presented an iterative flow updation algorithm to obtain user equilibrium and system optimal flows. Research since then has focused on the even more general asymmetric traffic assignment problem, and there have been no recent notable advances focusing specifically on symmetric TAP. A core contribution of this study is adaptation of Algorithm B for symmetric TAP by proving correctness of the flow shifts and ensuring that required optimality conditions are met.

The solution to static symmetric TAP is unique if the link performance functions are strictly monotone and the cost Jacobian is positive semidefinite. These conditions are satisfied by the cost functions we describe in the Cost Formulae section. This study uses our implementation of Algorithm B to solve this assignment problem. Algorithm B selects the shortest and longest paths to each destination node within a bush (acyclic graphs rooted at each zone), and shifts flow using Newton's method, iterating over all bushes till convergence. The source code is available at the SPARTA lab github repository~\cite{github} under the \textit{pd-word} branch.

\subsection{Network design problem background}

Network design problems (like many other bi-level problems) are intractable to solve exactly. Such methods (such as branch-and-bound or Benders decomposition) have been proposed, as in \citeN{LeBlanc}, \citeN{Chen}, \citeN{Drezner}, and \citeN{Long}.  However, the largest network tested in any of these studies has 40 nodes and 99 links.

Heuristic methods are standard for this class of problems.  Genetic algorithms, simulated annealing, and tabu search are examples of such methods that have been applied to traffic NDPs. For details on such methods, see the review papers by \citeN{Farahani2014} and \citeN{Iliopoulou2019}. They study the urban transportation NDPs and provide an overview of the types of problem variations as well as solution methods and some applications. 

Alternatively, the discrete NDP can be transformed into a single level problem. \citeN{Gao2005} transformed the problem into a nonlinear single level program utilizing support functions, which was solved by existing techniques. A third approach involves formulating the problem with equilibrium constraints, then using branch and bound or reduced gradient based methods~\cite{Lo2004,Lo2009,Szeto2006,Szeto2008,Szeto2010}. Other methods, such as Lagrangian relaxation and column generation, have also been used as exact solution methods for small NDP instances~\cite{Meng2001,Borndorfer2008}. As this approach is not viable for large instances, we focus on the non-exact methods, which trade off guaranteed optimality for tractability. 

In our experiments, we use a genetic algorithm as a solution heuristic. \citeN{katoch2021review} provide a recent review of GAs, with applications in transit NDPs, multimodal networks, facility layout, inventory and scheduling, etc.  Genetic algorithms have some key advantages useful for our application, namely, easy parallelization, usage of objective function information (as opposed to gradient/hessian information), and ability to quickly obtain quality solutions in practice. These advantages are especially suited to our application, given computational and storage size requirements for large networks and lack of accompanying first/second-order information.  The downsides of GA such as high computation times and efficient representation of solutions/operators can be reasonably remediated with computational resources and proper implementation practices.

\section{Model Formulation}
\label{sec:model_formulation}

Let $\mathbf{A}$ denote the set of rail links traversable by diesel-electric trains, and $\mathbf{N}$ denote the set of nodes (yards and stations). Each link $i \in \mathbf{A}$ currently has a flow-dependent usage cost $c_i (x_i)$ per unit flow.  The link can be electrified for a cost of $c^e_i$, changing the usage cost function to $c'_i (x_i)$ per unit flow. The usage costs for diesel and electrified links differ due to technological differences and fuel costs, and are lower for electrified links. The set of candidate links eligible for electrification is denoted by $\mathbf{A^E}$, a subset of $\mathbf{A}$. The flow on link $i$ is given by $x_i$ tons per day. The capacity of the link is denoted by $u_i$. Let $\Pi$ denote the set of paths in the network, and $h_\pi$ the flow on a particular path $\pi \in \Pi$. 

The set of nodes (denoting stations, yards, interchanges, etc.) is denoted by $\mathbf{N}$.
The demand between each origin $r$ and destination $s$ node is given by $d_{rs}$. The demand information between all origin-destination (OD) pairs is stored in the OD matrix $\mathbf{D}$. The total budget for upgrades is $B$. The decision variable $y_i$ for $i \in \textbf{A}$ equals 1 if link $i$ is chosen for upgrades (electrification) and 0 otherwise. The cost for switching the mode of operations from diesel to electric (and vice versa) at yard $u$ is $w_u$ (and infinitely high at other nodes without switching facilities).

\begin{figure}[h]
  \includegraphics[width=\linewidth]{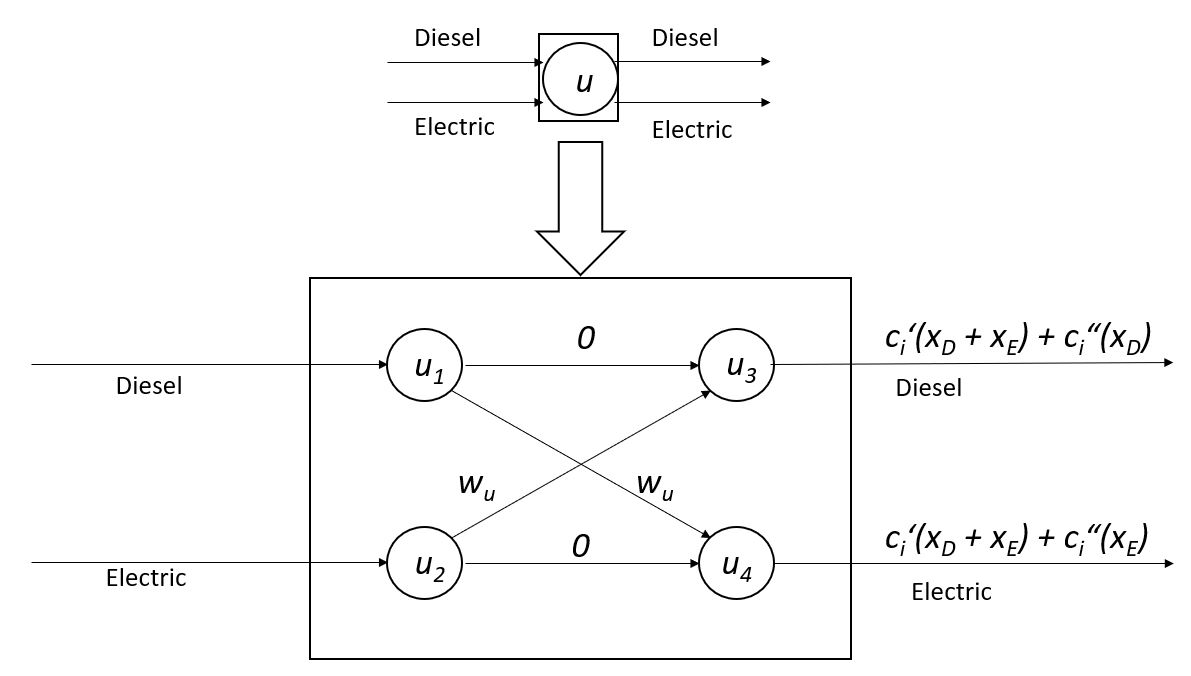}
  \caption{Network transformation for every node}
  \label{fig:net_transform}
\end{figure}

The practical interpretation of link electrification involves additional electric infrastructure on a link allowing electric trains to use the same physical link with different fuel costs, while maintaining the same total link capacity. Our model captures this using the network transformation shown in Figure \ref{fig:net_transform} with a node, an incoming link, and an outgoing link. Each physical link is split into two parallel sub-links, one for diesel flow ($x_{Di}$) and one for electric flow ($x_{Ei}$). The congestion cost $c_i'$ on each electrified link is not separable, depending on the sum of diesel and electric flows, thus capturing the physical link capacity constraint. The fuel (track resistance) costs $c_i''$ on these two sub-links are separable, depending only on the sub-links own flow (diesel or electric). The capacity of the electric sub-link is set infinitely high for links not chosen for electrification. Therefore, the total link cost of any link $i$ is $c_i(x_{Di} + x_{Ei}) = 2c_i'(x_{Di} + x_{Ei}) + c_i''(x_{Di}) + c_i''(x_{Ei})$. Each node $u$ in the network is expanded to incorporate switching costs $w_u$. This switching cost is set infinitely high at non-yard nodes to allow switching only at yards.

\noindent The RNDP is then formulated as follows:

\begin{align}
    \min \quad\ \ &\sum_{i \in \mathbf{A}}   c_i (x_i)  x_i \label{eqn:rail1}\\
    \nonumber\text{subject } &\text{to:}\\
    &\sum_{i \in \mathbf{A}} c^e_i y_i \leq B \label{eqn:rail2}\\
    &x_i =  x_{Di}+x_{Ei} \qquad \forall  i \in \mathbf{A}\\
    &y_i \in \{0,1\} \qquad \forall i \in \mathbf{A}\label{eqn:rail3}
\end{align}

\noindent where $y_i$ are found from the upper level problem and $x_i$ values are found by solving the following lower level assignment problem, in which $M$ is a sufficiently large constant:

\begin{align}
    \min_{\mathbf{x},\mathbf{h}} \quad\quad  &\oint_0^{\mathbf{x}} \mathbf{c} (\mathbf{s})  \cdot d\mathbf{s}\label{eqn:rail4}\\ 
    %\sum_{i \in \mathbf{A \cup A^E}}
    \nonumber\text{subject to:}&\\
    x_{i} &= \sum_{\pi \in \Pi : (r,s) \in \pi} h_\pi & \forall (r,s) \in \mathbf{A} \label{eqn:rail5}\\ %\cup A^E}
% SDB: A \cup A^E is redundant since A^E is a subset of A, right?    
    x_i & \leq M y_i   & \forall i \in \mathbf{A^E} \label{eqn:rail6}\\
    \sum_{\pi \in \Pi^{rs}} h_\pi &= d_{rs} & \forall(r,s) \in \mathbf{N}^2\label{eqn:rail7}\\
    h_\pi &\geq 0 & \forall \pi \in \Pi\label{eqn:rail8}
\end{align}

Equation~\ref{eqn:rail1} minimizes the total usage cost over electrified and non-electrified links. Equation~\ref{eqn:rail2} enforces the electrification budget constraint. Equation~\ref{eqn:rail3} indicates that $y_i$ is a binary indicator variable. Equation~\ref{eqn:rail4} is the modified Beckmann function for symmetric link interactions, described more in the next sub-section. Equation~\ref{eqn:rail5} states that the relation between link and path flows, and Equation~\ref{eqn:rail6} ensures that the electrified links with flow are the ones actually chosen for construction. Equation~\ref{eqn:rail7} ensures that the demand is met across all used paths for each OD pair, and Equation~\ref{eqn:rail8} states flow non-negativity.

We would like to make an important observation here. The total demand is considered constant for our case study, and we conduct preliminary sensitivity analysis to understand the effect of varying (induced) demand. The problem with elastic demand is equivalent to the TAP formulation via the Gartner transformation ~\cite{gartner1980optimal} if the model elasticity is fully known. This transformation creates direct links from every origin to destination and captures potential flow increase using the elasticity data on the new links. However, this requires model elasticity to be fully known, and the network size increases by several orders of magnitude. Therefore, despite being equivalent from a method standpoint, we have not incorporated elastic demand in our case study.

\section{Cost Formulae}

This section provides the cost calculation formulae for the link electrification costs, track resistance, and generalized link costs. Owing to a significant number of constants and rates within these formulae, they are omitted here for brevity, but can be found on the author's github repository~\cite{github2}.  In order to update cost values for inflation, this study uses industry appropriate producer price indices from the \citeN{bls}.

There are two main costs to be considered, roughly corresponding to the capital expenses and operational expenses. Electrification costs (or capital expenses) refer to all the costs associated with converting a link from standard operations to electric operations that we considered.  The generalized link costs (operational expenses) refer to the costs of traversing a link, and incorporate both time and fuel (whether diesel or electricity) costs based on the link's speed and resistance, calculated on a per-ton basis.

\subsection{Link Electrification Cost}

The link electrification cost does not depend on the link's traffic/tonnage (costs associated with re-routing during construction are not considered), but rather the link's location and existing characteristics. This cost can be divided into five components:

\begin{itemize}
    \item Overhead contact system (OCS) --- $c_{OCS}$
    \item Electrical substations --- $c_{sub}$
    \item Transmission lines --- $c_{trn}$
    \item Public works --- $c_{pub}$
    \item Signalling --- $c_{signal}$
\end{itemize}

The first four cost components are positively correlated with the roughness of the link's terrain (e.g. the electrification infrastructure is more expensive to construct in mountainous segments than in plains segments).  For each of those four categories, the cost calculations utilize two parameters: a low value for ideal conditions, and a high value for the worst possible condition.  In order to create an upgrade cost estimate for all links in the network, the difficulty of the terrain is scaled linearly based on the ratio, $\alpha$, between each link's actual length and its straight-line length.  $\alpha$ has a minimum value of one, which represents a straight track.  The cost of signaling depends on the complexity of replacing the link's existing signal systems, and is independent of terrain.

\noindent The cost of upgrading link $i$ to be electrified is then written as:

\begin{equation}
\begin{split}
c^e_i=l_i \left(\left(\frac{\alpha_i -\alpha_{min}}{\alpha_{max}-\alpha_{min}}\right)(c_{OCS,max} + c_{sub,max} + c_{transmission,max} + c_{pub,max}) \right. +\\
 \left. \left(1-\frac{\alpha_i -\alpha_{min}}{\alpha_{max}-\alpha_{min}}\right) (c_{OCS,min} + c_{sub,min} + c_{trn,min} + c_{pub,min})+c_{signal, i}\right)
\end{split}
\,,
\end{equation}
where $l_i$ is the length of link i, and $\alpha_{min}$ and $\alpha_{max}$ are the smallest and largest terrain difficulty values across the network, respectively. The parameters are derived from \citeN{schwarm1977}, \citeN{Kneschke1986}, and \citeN{gattuso2014tool}.

\subsection{Rail Link Delay Function}

To assign freight flows to the network, a link performance function is required, accounting for the rising cost of using a link as link congestion increases. This subsection details the formulation of the function as well as the process used to obtain the coefficients and constants. The final form of the equation is provided in equation \ref{eqn:GC1} at the end of this subsection.

Several studies have applied a rail link performance function similar in form to the \citeN{BPR} function used in road assignment~\cite{clarke1995}. \citeN{uddin2015} proposed a link delay function of the form:

\begin{equation}
t_i=t_i^0 \left(1+\left(\frac{x_i}{u_i}\right)^\beta\right)
\end{equation}

\noindent where $x_i$ is the flow on link $i$, $u_i$ is the link's capacity, and $t_i^0$ is the free-flow travel time. Uddin suggested a value of 4 for the parameter $\beta$. $u_i$ is a function of the link's class, as well as the number of tracks in the link and the frequency of sidings and switches.

The link's total generalized cost is a function of the travel time, crew and cargo costs, as well as the fuel cost. The fuel cost is the product of the link travel time, cost per unit energy for the fuel source used, and power level required for the link divided by the efficiency of the locomotive. For this analysis, an average train unit is the basis of analysis, and the generalized cost is based on the cost for that train unit to traverse a link divided by that train unit's cargo mass. In this way, the network flows can be assigned as tons of cargo.  The power level used on each link is calculated using the link's total resistance, as outlined in the next subsection. The generalized cost function can then be written as follows:

\begin{equation}
\label{eqn:GC}
c_i\ (or\ c'_i)= t_i \left(crewRate\ +\ cargoRate\ +\ \frac{P_i}{\eta} fuelCost \right)/(n_C m_{cargo})
\end{equation}
where $n_C m_{cargo}$ is the mass of the cargo hauled by a train unit used in the analysis.

\noindent In order to separate the congestion cost and the resistance cost, we make the assumption that the fuel consumption on each link is fixed, meaning the generalized cost function changes as follows:
\begin{equation}
\label{eqn:GC1}
c_i\ (or\ c'_i)= \left(t_i \left(crewRate\ +\ cargoRate\ \right) + t_0\left( \frac{P_i}{\eta} fuelCost \right)\right)/(n_C m_{cargo})
\end{equation}

%\noindent We assume that the power level does not change according to link delay.  This assumption can be relaxed in future work to represent, for instance, that a high link delay might allow the train to operate at a lower power level, or it could imply more brake usage.

\subsection{Track resistance}

%Assume an average train, in this paper refered to as a train unit.
In this paper, we assume that trains have $n_L$ locomotives and $n_c$ railcars. Each railcar has a tare weight of $m_c$ and a gross weight of $m_g=m_c+m_{cargo}$.  Each locomotive has a weight of $m_L$. A loaded train then has a mass of $m_T=n_L m_L+n_c m_g$.  
%Because the model is not able to include the cost of switching a consist from electrical locomotives to diesel-electric locomotives and vice versa, the average train in the simulation will have both sets of locomotives.  The locomotives that are not in use on a particular link do not contribute to the train's tractive effort, but they do contribute to the resistance as deadweight. 
Time and energy exertions associated with switching between diesel-electric and electric locomotives are handled separately, and discussed in the switching costs section.
%commented-out sentences related to prior analysis when switching costs couldn't be incorporated.  Added sentence to clarify that the track resistance related to switching movements are considered separately.
The resistance on the train is separated into bearing resistance, flange resistance, air resistance, grade resistance, curve resistance, and brake resistance, and inertial resistance. Each of these quantities is discussed and specified below.
\begin{description}
\item[Bearing resistance:]
Assuming a relatively new train (less than fifty years old), the bearing resistance on each railcar will be $(a+(\frac{b N_{ax}}{m_g}))m_g$, where $N_{ax}$ is the number of axles on the railcar. According to \citeN{AREMA2003}, when $m_g$ is measured in tons, $a=2.9 \frac{\mathrm N}{\mathrm {ton}}$ and $b=97.3N$. Therefore, the total bearing resistance is calculated as $\sum_{k=1}^{n_L+n_C} (a m_{g,k}+b N_{ax,k})$.

\item[Flange resistance:]
The bearing resistance on each component of the train is $m(B v)$, where $v$ is the speed of the train relative to the track, and $m$ is the mass of the railcar or locomotive~\cite{AREMA2003}. $B$ = 0.329 $\frac{\mathrm{N} \cdot \mathrm s}{\mathrm m \cdot \mathrm{ton}}$ for locomotives($B_L$), and 0.494 $\frac{\mathrm N \cdot \mathrm s}{\mathrm m \cdot \mathrm{ton}}$ for railcars($B_C$). The total flange resistance is $v(B_Ln_L+B_Cn_C)$. The flange resistance varies based on the track quality, so a factor $k_{f,i}$ can be applied to adjust the flange resistance for link $i$.

\item[Air resistance:]
The air resistance\cite{hay1982railroad} on each railcar is proportional to the square of the train's speed relative to the wind. For the purposes of this study, the wind speed is assumed to be zero, so that the train's speed relative to the track is equal to its speed relative to the air. The air resistance on each component of the train becomes $Kv^2$. When $v$ is measured in $\frac{\mathrm m}{\mathrm s}$, $K$ takes on a value of 1.56 $\frac{\mathrm N \cdot \mathrm{s}^2}{\mathrm{m}^2}$ for conventional equipment, or 2.06 $\frac{\mathrm{N} \cdot \mathrm{s}^2}{\mathrm{m}^2}$. The total air resistance is $\sum_{k=1}^{n_L+n_C} K_k v^2$.  A factor $k_{a,i}$ can be used to adjust the air resistance based on the average air density for link $i$.

\item[Grade resistance:]
The grade resistance is the resistance due to gravity. Unlike the other resistances discussed, grade resistance can be positive (upgrades) or negative (downgrades). The grade resistance for the train is given by $m_T g \sin(\theta)$, where $g$ is the acceleration due to gravity and $\theta$ is the angle of incline.  Rail inclines are small, allowing the small-angle approximation that $\sin(\theta)$ is approximately equal to the track grade $gr$.

\item[Curve resistance:]
Curve resistance arises from the force of the track on the wheels within a curve.  According to AREMA~\cite{AREMA2003}, curve resistance is approximately equivalent to a grade of 0.04\% per degree of curvature. That heuristic would put the curve resistance at $0.45836 m_T g \arcsin\left(\frac{30.48}{2r}\right)$, where $r$ is the link's radius of curvature in meters (an average radius of curvature is assigned to each link based on the link's $\alpha$ value).

\item[Brake resistance:]
The brake resistance is the force of the brakes applied to the train. This force is applied to maintain control of the train along down grades, and to stop the train.  Trains have three braking systems, the most powerful of which, the air brake, takes large amounts of time to engage or disengage.  %Sometimes, the air brake is left on along flat terrain or small up grades to prevent delays accruing from recharging the pressure in the train tube.  
Electric locomotives have regenerative braking systems that allow them to recover some braking energy.  The regenerative braking system is more reliable than a diesel-electric locomotive's rheostatic braking system, allowing electric locomotives to rely less so on the air brake.

The cost parameters for fuel efficiency and maintenance are derived from \citeN{Whitford1981}, \citeN{fritz}, and \citeN{Nektalova}. Because we model trains as flows, it is beyond our scope to determine the actual brake usage on each link.  We assume that on level terrain, links with an average positive grade, or links with a negative grade below a certain threshold, the incidental brake usage will be equivalent to the resistance of a $0.1\%$ grade.  A threshold is determined using the grade that would cause the train to exceed its desired speed along the link when utilizing the minimum throttle level.  Beyond that threshold, the brake force is set to the level that will allow the train to reach its desired speed along the link while utilizing the minimum throttle position.  The throttle is not set to zero because the grade utilized is only the average along the link, and even when a train is going downhill for a significant distance, the motors can be operating to prevent the railcars from bunching together.

\item[Inertial resistance:]
Inertial resistance is the positive or negative impedance from the train changing its velocity.  Trains must use energy to accelerate, and much of that energy is not regained when slowing down due to friction or the need to slow down faster than would otherwise be necessary. The inertial resistance is $m_T a$.  This study assumes that the positive and negative inertial resistances will cancel each other out.  In reality, more energy is used in accelerating a train than the energy saved as a train decelerates, so this is one category where the total energy consumption is understated and future analysis could improve upon.

\item[Total resistance:]
Combining the preceding formulae yields the total average resistance along a link $i$ as:

\begin{equation}
    \label{eqn:resistance}
\begin{split}
R_i = \sum_{k=1}^{n_L+n_C} &(a m_{g,k} + b N_{ax,k} + K_k k_{a,i} v_i^2) + v_i(B_Ln_L + B_Cn_C)k_{f,i} \\ &+ m_T g (gr_i) + 0.4536m_T g \arcsin\left(\frac{15.24}{r_i}\right) + R_{brake, i} + m_T a
\end{split}
\end{equation}

where $k$ refers to each rail vehicle in the train unit.\\
\end{description}

\noindent \textbf{Power level and speed}: The power level used along the link, $P_i$ is a function of the link's resistance and the train's speed along the link:
\begin{equation}
    \label{eqn:power}
    P_i = R_i v_i
\end{equation}
The train unit has a discrete number of possible power levels, which it uses to approach the desired speed along each link. Equations~\eqref{eqn:resistance} and ~\eqref{eqn:power} are solved iteratively to determine each link's associated power level, $P_i$, and base travel time, $t_i^0=l_i/v_i$.  Those values are used in equation~\eqref{eqn:GC1}.

\subsection{Switching costs}
For many O-D pairs, the lowest-cost route for the shipper may involve switching between electric and diesel-electric operations en-route.  This might occur when a small portion of the overall network is electrified and an electrified path does not exist between on O-D pair, or the electrified route is far enough away from the direct path that the savings from electric operations are not worth the cost in added time.  If the network does not allow any switching, electric trains would only be assigned between O-D pairs that have fully electrified, mostly direct, paths between them, which is not realistic.  The network allows flows to switch from electric to diesel-electric operations (and vice versa) at certain nodes in order to simulate more realistic routing.

461 nodes in the rail network have facilities adequate to allow for switching.  At those nodes, links between the electric and diesel-electric links represent the cargo time, crew time, and energy costs associated with switching locomotives.  The cost assumes an hour-and-a-half per switch, and that six employees will be involved (including the train operators and the yard workers).  Whenever the path switches between operations, each train accrues a cost of roughly $\$3800$. This cost represents in-part the cost of yard electrification, i.e., additional switching using locomotive caller switchers.

%_________________________________________________________________________

\section{Subproblem solution}
\label{app:algBconvergence}

Most algorithms for TAP are designed for the separable cost case.  In this sub-section, we describe the changes that need to be made for our cost structure, and derive a new flow shift formula for this case.  Our solution method is based on Algorithm B, briefly described in the literature review.  The main points of this algorithm are disaggregating link flows by origin, where the sets of used links form acyclic subnetworks; shifting flows from longer-cost segments to shortest-cost segments within each subnetwork; and updating these subnetworks by removing unnecessary links and adding new ones which provide shorter paths.

Of these three components, the original derivation of only the second component (flow shift) relies on separable link costs, and must be re-examined in light of interactions.  We first show that link interactions in our formulation satisfy the symmetry condition needed for the formulation~\eqref{eqn:rail4}--\eqref{eqn:rail8} to hold, and then re-derive the flow shift formula for Algorithm B with interactions.

As discussed in the Cost Formulae section, all flow continuing on the same fuel type has zero switching cost, whereas all flow changing fuel type has pre-defined switching costs. Note that only the parallel diesel-electric link pair has cost function dependent on the sum of the flows (see Figure \ref{fig:net_transform}), and $c'(x)$ is identical for both links. The Jacobian takes the form:

\begin{equation*}
J_{\mathbf{x}} = 
\begin{pmatrix}
\frac{\partial(c_1(\mathbf{x}))}{\partial x_1} & \frac{\partial(c_2(\mathbf{x}))}{\partial x_1} & \cdots & \frac{\partial(c_n(\mathbf{x}))}{\partial x_1} \\
\frac{\partial(c_1(\mathbf{x}))}{\partial x_2} & \frac{\partial(c_2(\mathbf{x}))}{\partial x_2} & \cdots & \frac{\partial(c_n(\mathbf{x}))}{\partial x_2} \\
\vdots  & \vdots  & \ddots & \vdots  \\
\frac{\partial(c_1(\mathbf{x}))}{\partial x_n} & \frac{\partial(c_2(\mathbf{x}))}{\partial x_n} & \cdots & \frac{\partial(c_n(\mathbf{x}))}{\partial x_n} \\ 
\end{pmatrix}
\end{equation*}

The off-diagonal elements are zero whenever the two links do not interact. In our case, the only non-zero off diagonal elements are for parallel diesel-electric pairs. These elements are calculated using the chain rule:

\begin{equation}
\frac{\partial c'(x_D + x_E) + c"(x_D)}{\partial x_E} = \frac{\partial c'(x_D + x_E)}{\partial(x_D+x_E)}.\frac{\partial(x_D+x_E)}{\partial x_E} = \frac{\partial c'(x_D + x_E)}{\partial(x_D+x_E)}
\end{equation}
\begin{equation}
\frac{\partial c'(x_D + x_E) + c"(x_E)}{\partial x_D} = \frac{\partial c'(x_D + x_E)}{\partial(x_D+x_E)}.\frac{\partial(x_D+x_E)}{\partial x_D} = \frac{\partial c'(x_D + x_E)}{\partial(x_D+x_E)}
\end{equation}

As can be seen, these two derivatives are identical, proving that the Jacobian is a symmetric matrix, and establishing the validity of the formulation~\eqref{eqn:rail4}--\eqref{eqn:rail8}.

We next show that the Algorithm B flow shift procedure (equalizing cost between routes) is valid even in the presence of link interactions.  This result is more general, and can be applied to \emph{any} instance of symmetric costs depending on at most two links, not just the ones we adopt in our study.  We separate out the terms involving interactions as follows: $c_i = f_1(x_i, x_j) + f_2(x_i)$ and $c_j = g_1(x_i, x_j) + g_2(x_j)$, with $\frac{\partial f_1(x_i, x_j)}{\partial x_j} = \frac{\partial g_1(x_i, x_j)}{\partial x_i}$. For our rail electrification formulation, $f_1(x_i, x_j) = g_1(x_i, x_j) = c_i'(x_D,x_E)$, resulting in a simpler problem.

In this flow-shift operation, we have identified two paths (a lower-cost path $\pi_L$ and a higher-cost $\pi_U$), and wish to shift flow between them to minimize the objective.  Recall that the line integral formulation can now be written as:
\begin{equation}
    F(\mathbf{x}) = \oint_{\mathbf{0}}^{\mathbf{x}} \mathbf{c}(\mathbf{s}) \cdot d\mathbf{s}
\end{equation}
Partition the arc set $\mathbf{A}$ into three subsets $A_1$, $A_2$, and $A_3$. Links in set $A_3$ have separable cost, depending only on their own flow. The travel time on the links in sets $A_1$ and $A_2$ depend on the flows on two links: the link itself, and exactly one link in the other set. Mathematically, there is a bijection between sets $A_1$ and $A_2$, with the notation $i(j)$ and $j(i)$ used to denote the counterparts of links in the other set. The notation for these cost functions lists the links own flow first, and its counterpart's flow second, that is, $c_i(x_i, x_{j(i)})$ and $c_j(x_j, x_{i(j)})$. Let $p = |A_1| = |A_2|$ and $m = |A_1| + |A_2| + |A_3|$. The index $a$ will be used to denote a generic link, if it doesn't matter which set it's from.

% For a given vector of link flows $\mathbf{x}$, the objective function is given by
% \[
% F(\mathbf{x}) = \int_{\mathbf{0}}^{\mathbf{x}} \mathbf{t}(\mathbf{s}) \cdot d\mathbf{s}
% \,.
% \]
As the line integral is path-independent, we choose the following integration path: 
\[
(0, 0, 0, \ldots,0) \rightarrow
(x_1, 0, 0, \ldots, 0) \rightarrow
(x_1, x_2, 0, \ldots, 0) \rightarrow
\cdots \rightarrow
(x_1, x_2, \ldots, x_m)
\,.
\]
The line integral then decomposes into a sum of ordinary integrals:
\begin{align*}
F(\mathbf{x})
=&
\sum_{i=1}^p \int_{(x_1, \ldots, x_{i-1}, 0, 0, \ldots, 0)}^{(x_1, \ldots, x_{i-1}, x_i, 0, \ldots, 0)} c_i(\mathbf{s}) \cdot d\mathbf{s}
+
\sum_{j=p+1}^{2p} \int_{(x_1, \ldots, x_{j-1}, 0, 0, \ldots, 0)}^{(x_1, \ldots, x_{j-1}, x_j, 0, \ldots, 0)} c_j(\mathbf{s}) \cdot d\mathbf{s} \\
&+ \sum_{k=2p+1}^{m} \int_{(x_1, \ldots, x_{k-1}, 0, 0, \ldots, 0)}^{(x_1, \ldots, x_{k-1}, x_k, 0, \ldots, 0)} c_k(\mathbf{s}) \cdot d\mathbf{s} \\
=&
\sum_{i=1}^p \int_0^{x_i} c_i(s,0)~ds
+
\sum_{j=p+1}^{2p} \int_0^{x_j} c_j(s,x_{i(j)})~ds
+
\sum_{k=2p+1}^m \int_0^{x_k} c_k(s)~ds 
\,.
\end{align*} 

Let $I_a$ be $+1$ if link $a$ is on $\pi_L$, $-1$ if $a$ is on $\pi_U$, and $0$ otherwise. Then for a flow shift $\Delta x$ from the longest to the shortest path, the change on each link's flow is given by $\Delta x_a = I_a \Delta x$ and $I_a = \partial (\Delta x_a)/\partial (\Delta x)$.
(This latter equation will be used when differentiating via the chain rule below). The objective in terms of a flow shift of size $\Delta x$ is now
\begin{multline*}
F(\Delta x)
=
\sum_{i=1}^p \int_0^{x_i + \Delta x_i} c_i(s,0)~ds
+
\sum_{j=p+1}^{2p} \int_0^{x_j + \Delta x_j} c_j(s,x_{i(j)} + \Delta x_{i(j)})~ds
+
\sum_{k=2p+1}^m \int_0^{x_k} c_k(s)~ds 
\,,
\end{multline*}
and when we differentiate with respect to $\Delta x$, we obtain the following (the second term splits into two terms from the Leibniz rule):
\begin{multline*}
\frac{dF}{d (\Delta x)} 
=
\sum_{i=1}^p c_i(x_i + \Delta x_i,0) I_i
+
\sum_{j=p+1}^{2p} c_j(x_j + \Delta x_j, x_{i(j)} + \Delta x_{i(j)}) I_j
\\+
\sum_{j=p+1}^{2p} \int_0^{x_j + \Delta x_j}
    \frac{\partial c_j}{\partial x_{i(j)}} (s, x_{i(j)} + \Delta x_{i(j)}) I_{i(j)}~ds
+
\sum_{k=2p+1}^m c_k(x_k + \Delta x_k) I_k
\,.
\end{multline*}
We will show that $dF/d(\Delta x) = \sum_{a\in(\mathbf{A}\cup \mathbf{A^E})} c_a(\mathbf{x} + \Delta \mathbf{x}) I_a$.  The interpretation is that the derivative vanishes if $\Delta x$ equalizes the cost difference on the longest and shortest segments, thus minimizing the Beckmann function. The second and fourth terms in the sum are exactly what we need for the links in $A_2$ and $A_3$.
We need to show \begin{multline*}
\sum_{i=1}^p c_i(x_i + \Delta x_i,0) I_i 
+
\sum_{j=p+1}^{2p} \int_0^{x_j + \Delta x_j}
    \frac{\partial c_j}{\partial x_{i(j)}} (s, x_{i(j)} + \Delta x_{i(j)}) I_{i(j)}~ds
\\=
\sum_{i=1}^p c_i(x_i + \Delta x_i, x_{j(i)} + \Delta x_{j(i)}) I_i
\,.
\end{multline*}
Using the symmetry condition $\frac{\partial c_j}{\partial x_{i(j)}} = \frac{\partial c_i}{\partial x_{j(i)}}$, we have
\begin{align*}
\sum_{i=1}^p c_i(x_i + \Delta x_i,0) I_i 
+
\sum_{j=p+1}^{2p} \int_0^{x_j + \Delta x_j}
    \frac{\partial c_i}{\partial x_{j(i)}} (s, x_{i(j)} + \Delta x_{i(j)}) I_{i(j)}~ds
\,.
\end{align*}

Using the fundamental theorem of calculus and one-to-one correspondence of terms in the two summations, it is clear that it equals $\sum_{i=1}^p c_i(x_i + \Delta x_i, x_{j(i)} + \Delta x_{j(i)}) I_i$, which is the required form. Therefore, $dF/d(\Delta x) = \sum_a c_a(\mathbf{x} + \Delta \mathbf{x}) I_a$ and Algorithm B flow shifts are valid in our setting.
Most implementations of Algorithm B (including ours) use a Newton method to equalize the costs on these links.  The denominator in Newton's method (the flow shift scaling factor) can also be adjusted to reflect interactions\cite{stap}, but even without this change, existing implementations still converge to the equilibrium.

\section{Experiment design}
\label{sec:datasources}

After obtaining solutions and validating the solution method for our formulation, the second set of experiments involves policy testing and sensitivity analysis.  These experiments vary parameters such as total budget, demand data, electrification costs, electricity costs, crew/cargo costs, policy changes in the form of monetary incentives. The base electrification budget is \$30 billion, roughly equivalent to electrifying $65,000$ kilometers of track. This is varied by up to $\pm 20\%$, or \$24--36 billion. The analysis involves studying the return on investment (ROI) in the form of reduced costs, incentivizing policymakers and private stakeholders to upgrade infrastructure.

Freight Analysis Framework, version 4 (FAF4) provides demand data from 2010 and 2020, as well as forecasts for 2030 and 2040~\cite{faf4}. In addition, lowering the usage cost of rail freight allows for modal shift from trucking to rail. This potential demand variation is considered by increasing the demand data by up to 25\% from base values. Lastly, electrification costs and crew/cargo costs are also tested at increased values (+25\%) to gauge effects on the electrification solution. Note that under the current experiment formulation, adjusting the electrification costs is equivalent to adjusting the budget for electrification.

\begin{figure}
     \centering
     \begin{subfigure}[b]{0.9\textwidth}
         \centering
         \includegraphics[width=\textwidth]{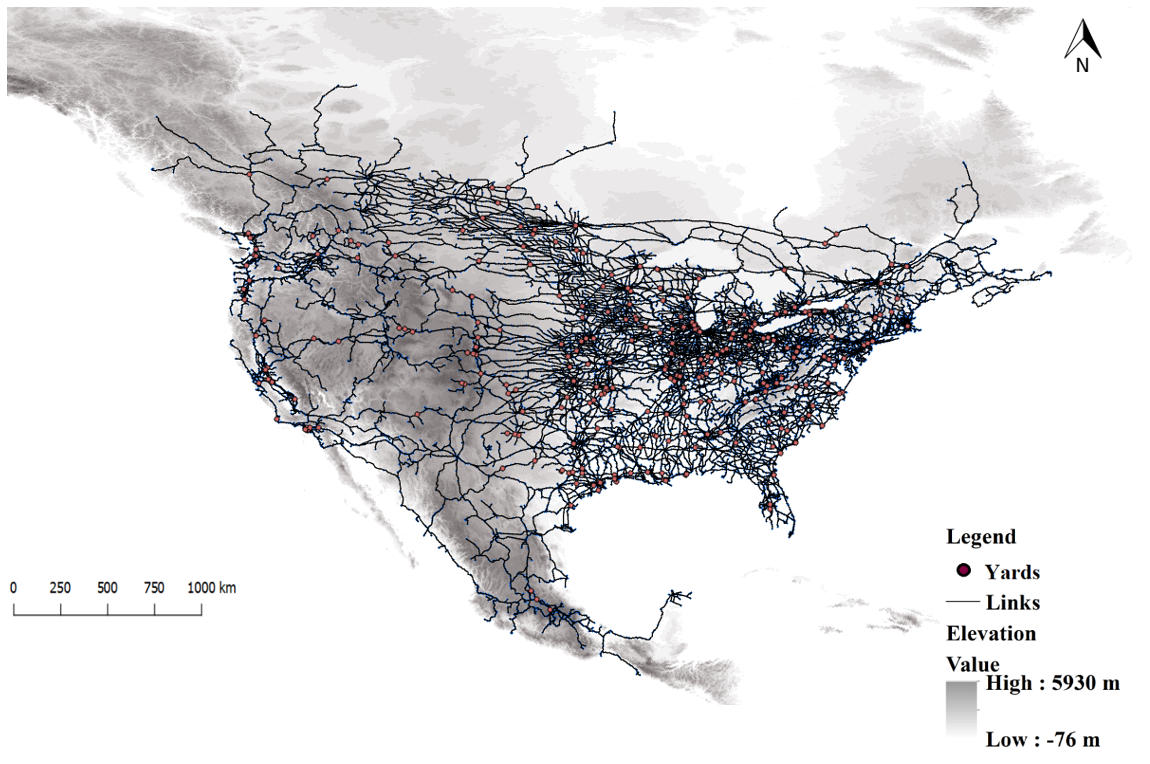}
         \caption{Base case results with shortest path network}
         \label{fig:basemap}
     \end{subfigure}
     \hfill
     \begin{subfigure}[b]{0.9\textwidth}
         \centering
         \includegraphics[width=\textwidth]{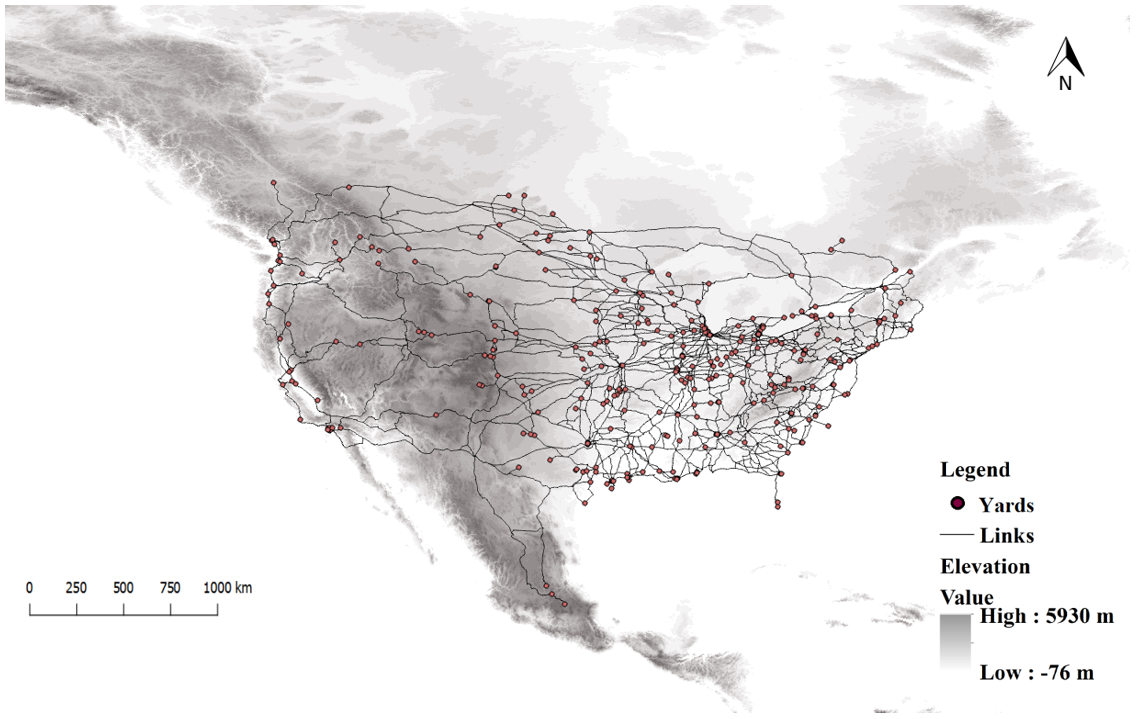}
         \caption{Candidate electrification corridor visualization between yards}
         \label{fig:candidate_corridors}
     \end{subfigure}
        \caption{Base network visualizations}
        \label{fig:base_net_viz}
\end{figure}

The North American railroad network has been extracted from the statewide analysis module~\cite{SAM} provided by the Texas Department of Transportation, originally based on the CTA rail network developed at Oak Ridge National Laboratory~\cite{CTA}. The dataset includes 35,424 links, and has geographical information (length, latitude, longitude, grade category), ownership information, and other auxiliary information. There are 28,289 nodes connected by these links, denoting stations, yards and interchanges. Elevation data were obtained by overlaying the network on the North American elevation grid~\cite{USGS} using ArcGIS (Figure \ref{fig:basemap}). The demand data (in tons) has been obtained from FAF4~\cite{faf4}. We would like to note that directional running (like implemented by Union Pacific in Texas) can result in increased flow efficiency, and can be modeled using increased link density. Barring detailed information on directional running, we have not included it in our model.

% \begin{figure}
% \includegraphics[scale=0.57]{Figures/2b_final.PNG}
% \centering
% \caption{Complete North American rail network representation with elevation data}
% \label{fig:basemap}
% \centering
% \end{figure}

As specified, operational mode changes can only occur at yards. Therefore, corridor electrification is more sensible than individual link electrification. Candidate corridors were obtained by connecting each yard to the nearest neighboring yards. Specifically, the following method was employed to find the set of candidate corridors and results can be seen in Figure \ref{fig:candidate_corridors}. This new shortest corridor network has tracks totaling $170,000$ kms of track available for electrification, down from over $305,000$ kms in the full network. The switching costs were obtained using the following assumptions for the locomotives: 1.5 hours to switch, throttle position 1 for diesel-electrics, 10\% power for electrics, and 6 crew-equivalents manpower.

\begin{enumerate}
    \item Calculate all pair shortest paths (APSP) from all yards
    \item Calculate APSP from all yards, but prune search branch after reaching a yard.
    \item Compare the two sets and keep the paths common to both sets
\end{enumerate}

% \begin{figure}[h]
%   \includegraphics[width=\linewidth]{Figures/2a_final.PNG}
%   \caption{Complete North American rail network representation with elevation data}
%   \label{fig:candidate_corridors}
% \end{figure}

%________________________________________________________________________

\section{Results and Sensitivity analysis}
\label{sec:results}

\begin{figure}
     \centering
     \begin{subfigure}[b]{0.9\textwidth}
         \centering
         \includegraphics[width=\textwidth]{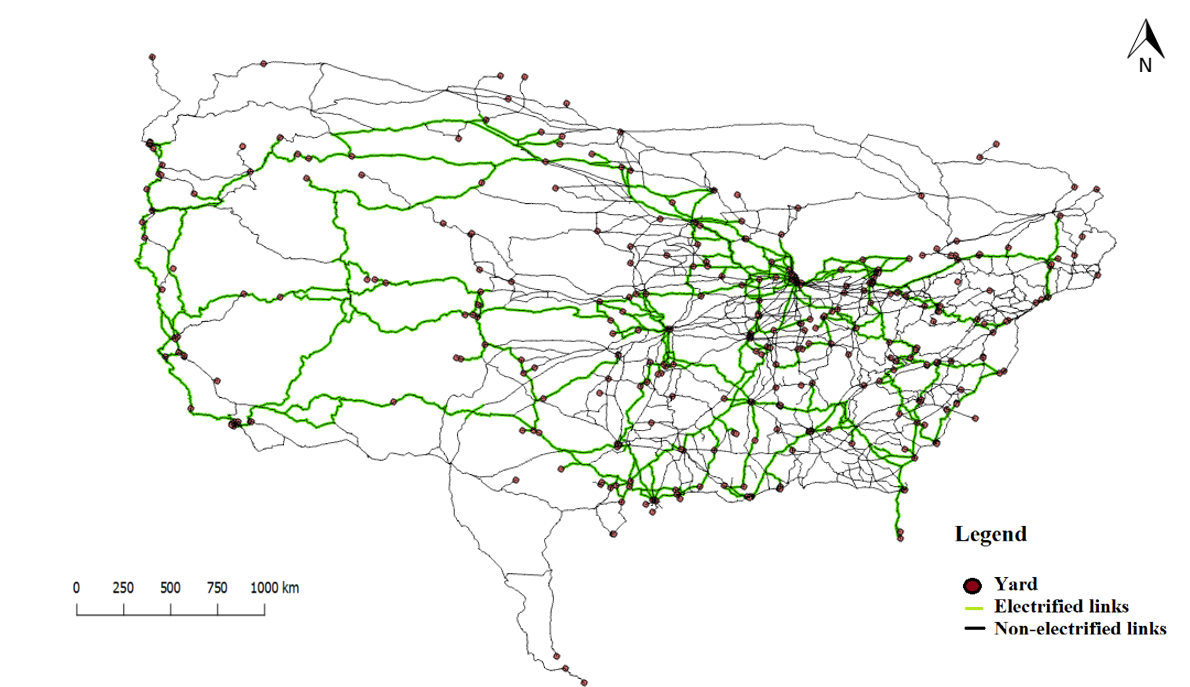}
         \caption{Base case results with shortest path network}
         \label{fig:base_case_shortest_paths}
     \end{subfigure}
     \hfill
     \begin{subfigure}[b]{0.9\textwidth}
         \centering
         \includegraphics[width=\textwidth]{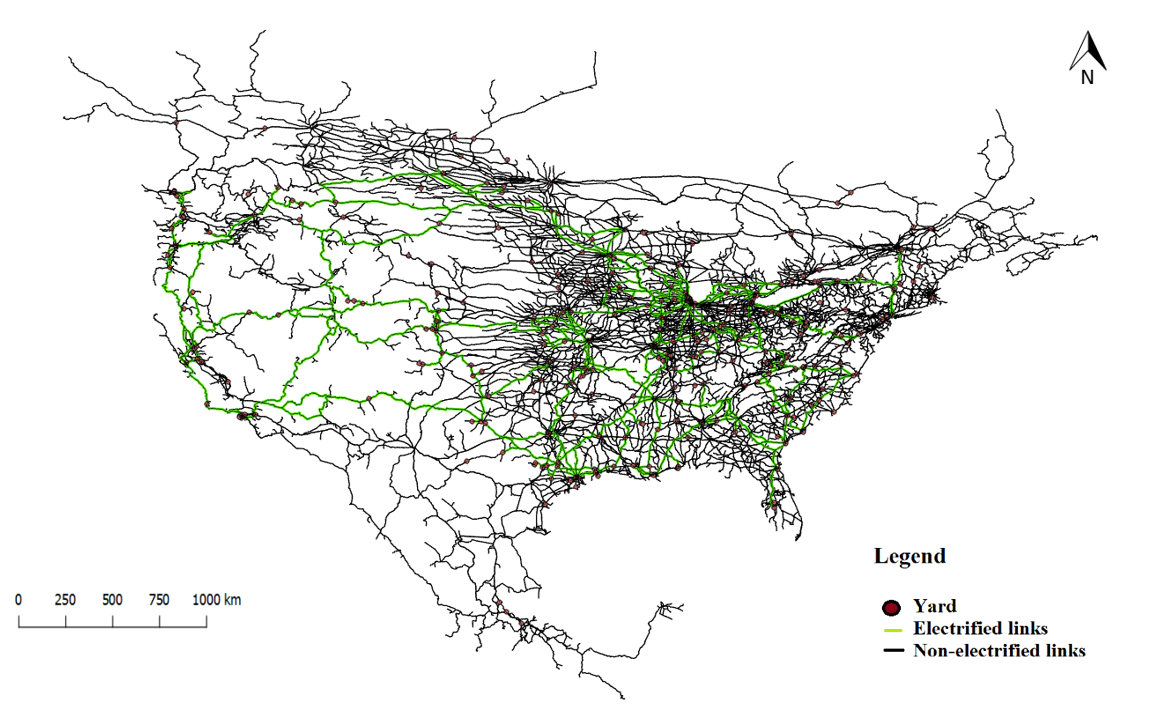}
         \caption{Base case results with full network}
         \label{fig:result_base_case_full_net}
     \end{subfigure}
        \caption{Base case visualizations}
        \label{fig:result_base_case}
\end{figure}

\begin{figure}[h]
  \includegraphics[width=\linewidth]{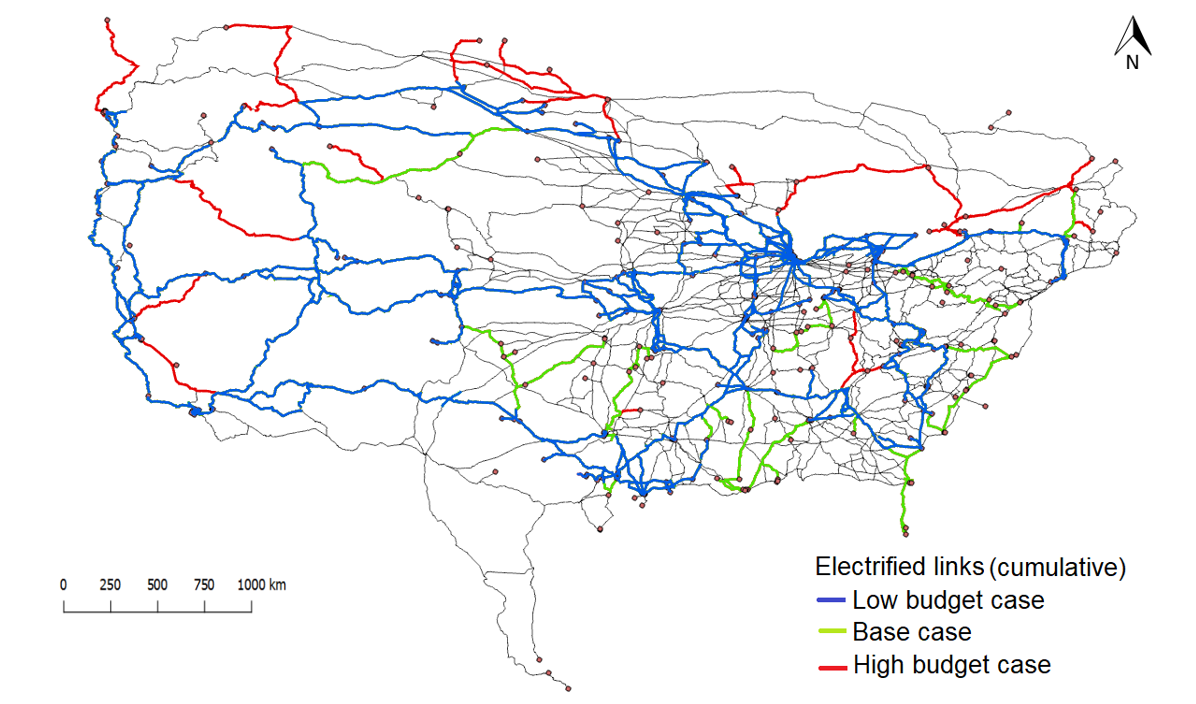}
  \caption{Results for varying electrification budget}
  \label{fig:results_varying_budgets}
\end{figure}

\begin{figure}[h]
  \includegraphics[width=\linewidth]{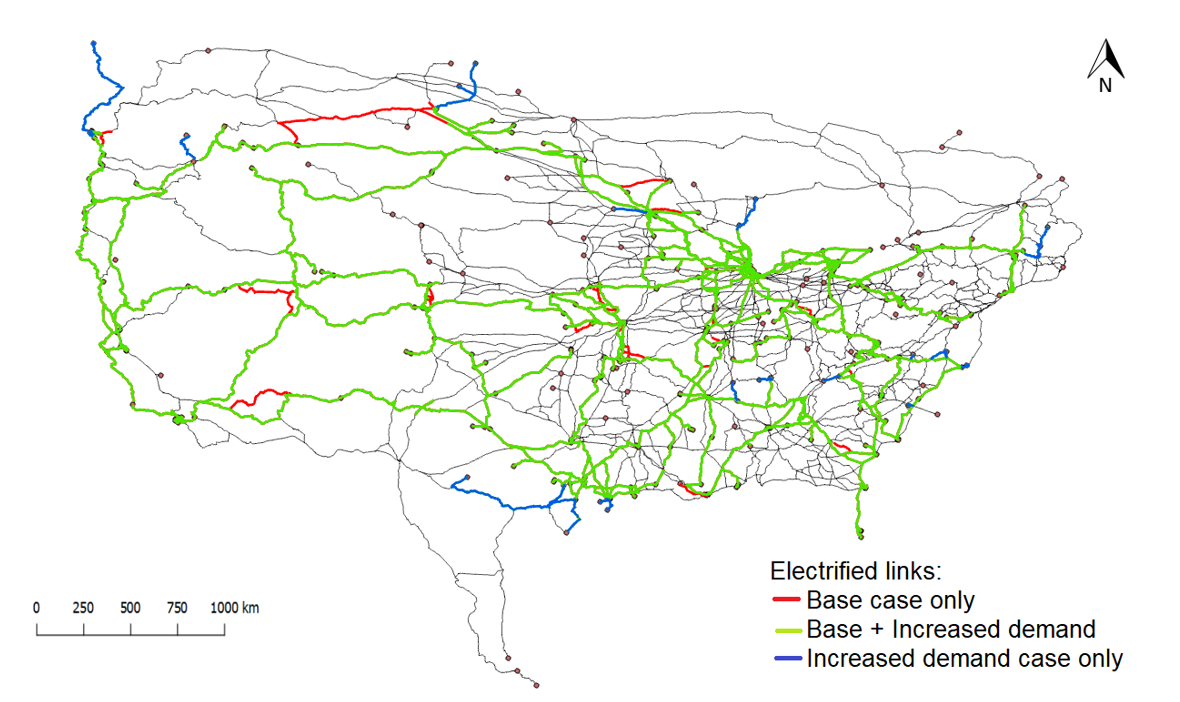}
  \caption{Results for increased demand case}
  \label{fig:results_increased_demand}
\end{figure}

\begin{figure}[h]
  \includegraphics[width=\linewidth]{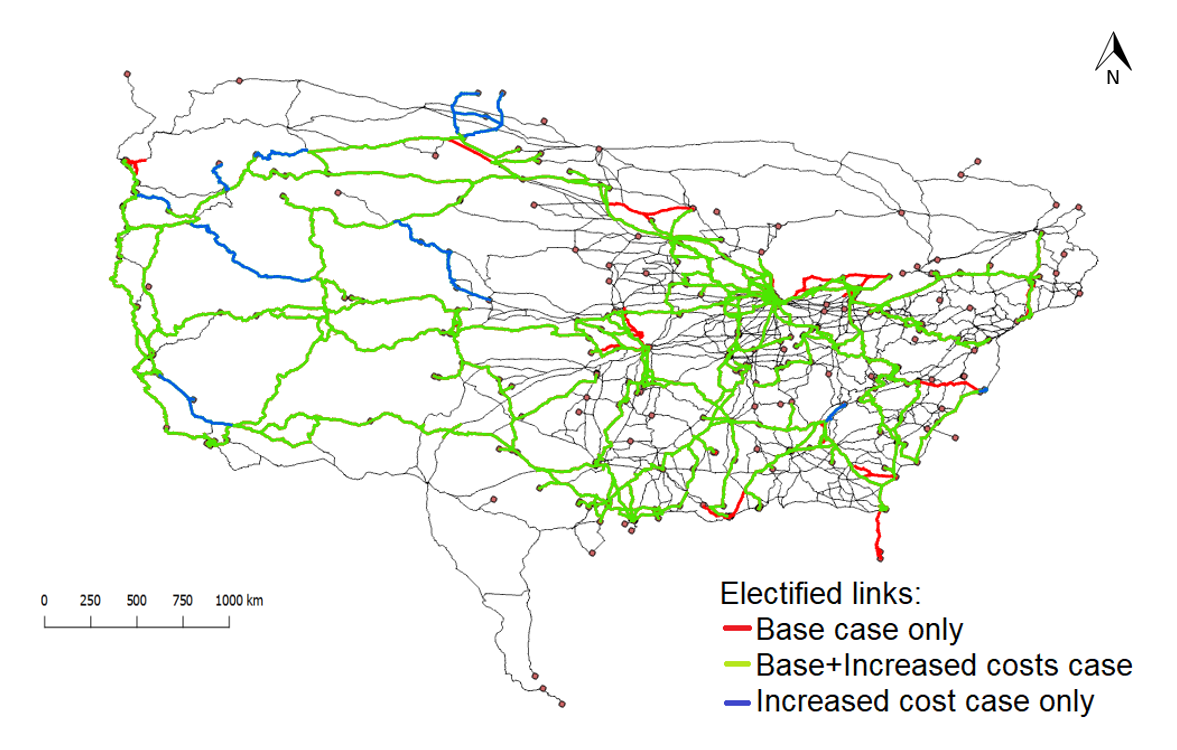}
  \caption{Results for increased operations cost case}
  \label{fig:results_increased_costs}
\end{figure}

All result visualizations are shown in Figures \ref{fig:result_base_case}--\ref{fig:results_increased_costs}. The low budget, base case (medium budget), and high budget cases lead to electrification of about 49,000 (30,000), 61,000 (36,000), and 70,500 (42,000) kms (mi), respectively. Only the base case results are shown with the full network backdrop, all other results are shown on the shortest path corridor network for better visibility. Figures \ref{fig:results_varying_budgets} - \ref{fig:results_increased_costs} show the results as the overlap of different scenarios. Figure \ref{fig:results_varying_budgets} is a cumulative plot, where the medium budget (base case) and high budget results are shown as increments over the low budget case. Figures \ref{fig:results_increased_demand} and \ref{fig:results_increased_costs} pivot off the base case, therefore, we highlight the overlapping links as well as the differences.

%Low_budget: 29796 miles
%Regular_budget:35984 miles 
%High_budget: 42063 miles

The first observation is that the GA retained all corridors selected for the low budget case in the base case and high budget case. However, this changes for the increased demand and increased costs scenarios, given that the network OD matrix and costs change in the two scenarios. This allows us to identify the most impactful corridors and stations, such as all three of the major transcontinental routes (LA-Chicago, Oakland-Chicago, and Seattle-Chicago), which the algorithm chooses for electrification in each scenario. There is a wide variety of corridors selected across the entire network. The algorithm selects quite a few mountainous routes where electrification provides the most benefits per train, but the algorithm also appears to select for links that provide connectivity throughout the network. This trend is best illustrated in figure 6, which seems to show that the increased demand case shifts the selection from the more mountainous west to the more populous east and gulf coasts. With higher demand, links that provide smaller benefits per train are more likely to provide greater benefits overall.

Examining the base case in further detail, the algorithm selects 13.2\% of the network (by line-miles) for electrification, and predicts that 15.5\% of final flow (by tonnage) will use electrical traction. This reflects that busier corridors benefit more from electrification. This also implies that most of the traffic along electrified links would be using electric rather diesel-electric locomotives.

The increased demand scenario increases electrified corridor connectivity to several new yards/stations. The increased cost scenario reduces selection from east and central US, and chooses corridors from western US. We hypothesize this selection occurs due to the higher elevations and grades for these tracks seeing benefit from electrification.  The current formulation does not incorporate any regional variations in the wholesale cost of electricity or diesel.

Each sub-problem took 5 minutes to converge, and each GA instance took ~36 hours to finish on a Linux PC with i7 Processor (2.6 GHz) and 16 GB RAM. Given the long-term planning horizon, runtimes are not crucial, but still useful to scalability of the heuristic for a large network.

%_________________________________________________________________________
\section{Conclusions and future directions}
\label{sec:conclusions}
This study presents a novel way of looking at the rail network electrification problem using connections to the road traffic assignment problem. We formulate the problem as a bi-level network design problem, where the lower level problem is a symmetric traffic assignment problem, assigning goods flow instead of traffic. We then show the correctness of Algorithm B flow shift formula for the symmetric traffic assignment problem, using it to solve this sub-problem. The costs and network parameters incorporate electrification costs; fuel, locomotive, and operational costs; and train resistance costs.

The North American railroad network is used as a test network to demonstrate our method, show scalability for large networks, and draw insights. We observe and note some corridors chosen in all different testing scenarios (varying budgets, opex, and demand), as well as provide policy insights for planning. The key observations are as follows:

\begin{itemize}
    \item While there seems to be some evidence of corridors in more mountainous terrain producing better cost savings, overall there is a wide variety of corridors selected across the entire network.  This implies that providing broad connectivity might be more important than electrifying the corridors with the best savings on a per-train basis.  
    \item The results for the increased demand case, shown in figure \ref{fig:results_increased_demand}, indicate a shift in chosen corridors from the more mountainous west towards the more populous east and gulf coasts.  This is probably a result of the increased demand causing those corridors to generate higher savings, even though they have fewer savings per train.  This is also evidence that the highest priority corridors might be routes through rough terrain operating at or near capacity.
    \item Sensitivity to demand means long-term trends in the trucking industry could have a substantial influence on the highest priority corridors for electrification.
\end{itemize}

The main limitations of this work are: 1) static power levels (once calculated); 2) static brake usage assumption on links; 3) static electricity and diesel costs; and 4) lack of incorporation of track ownership restrictions on path selection. The first three limitations are relatively straightforward to tackle, although they increase problem complexity and computation time significantly by introducing another layer of iterative calculations. 

An additional complication is partly based on the data available: because this analysis uses large zones with few centroids, the simulations might overestimate the benefit of electrifying well-traveled corridors that connect multiple zones. Studying how the selection changes based on varying switching costs might show whether the results are skewed by simulating a higher amount of connectivity than would actually exist.  This problem could partially be ameliorated by applying a cost to electric trains reaching centroid connectors to reflect that some electric trains might require additional drayage over diesel-electric trains. Properly calibrating such a cost would remain a difficult problem.

\label{sec:future directions}
The methods formulated in this study can be adapted to consider social benefits from changes in emissions, formulating the problem as one of social benefit maximization rather than private cost minimization. The composition of the power grid affects the emissions an electric locomotive causes, meaning the marginal social benefits of electrification and the marginal private cost savings can vary substantially from link to link.  In some parts of the power grid where most of the electricity is provided by coal, an electric train might even produce more emissions than an equivalent diesel-electric train~\cite{walthall}.  Understanding social benefits is important because the high variability of private return-on-investment from electrification might necessitate public subsidies before capital construction becomes possible. Prior studies have formulated social benefit maximization problems as a single-level problems, thus reducing computational complexity.

\section{Data Availability Statement}

All data, models, and code generated or used during the study appear in the submitted article.

\section{Acknowledgements}
The authors would like to thank William Alexander and Karthik Velayutham for their insights and help with the TAP code implementation.  This research was partially supported by the National Science Foundation (CMMI-1562291, CMMI-1826320).

\bibliography{ascexmpl-new}

\begin{thebibliography}{}

\bibitem[\protect\citeauthoryear{}{{Alliance Transportation Group}}{2013}]{SAM}
{Alliance Transportation Group} (2013).
\newblock ``{Statewide Analysis Model - Third Version (SAM-V3)}.

\bibitem[\protect\citeauthoryear{}{Bar-Gera}{2010}]{bar2010}
Bar-Gera, H. (2010).
\newblock ``Traffic assignment by paired alternative segments.''\ {\em
  Transportation Research Part B: Methodological}, 44(8-9), 1022--1046.

\bibitem[\protect\citeauthoryear{}{Beckmann et~al.\@}{1956}]{Beckmann}
Beckmann, M., McGuire, C.~B., and Winsten, C.~B. (1956).
\newblock ``Studies in the economics of transportation.''\ {\em Report No.
  RM-1488-PR}, RAND corporation.

\bibitem[\protect\citeauthoryear{}{Bornd{\"o}rfer
  et~al.\@}{2008}]{Borndorfer2008}
Bornd{\"o}rfer, R., Gr{\"o}tschel, M., and Pfetsch, M.~E. (2008).
\newblock ``Models for line planning in public transport.''\ {\em
  Computer-aided systems in public transport}, Springer,  363--378.

\bibitem[\protect\citeauthoryear{}{Boyles et~al.\@}{2020a}]{github}
Boyles, S., Alexander, W., Patel, R., Velayutham, K., and Thakkar, R. (2020a).
\newblock ``tap-b implementation.''\ {\em
  https://github.com/spartalab/tap-b/wrap}\ Accessed: 2020-05-20.

\bibitem[\protect\citeauthoryear{}{Boyles et~al.\@}{2020b}]{blubook_vol1_v085}
Boyles, S.~D., Lownes, N.~E., and Unnikrishnan, A. (2020b).
\newblock {\em Transportation Network Analysis}, Vol.~1.
\newblock 0.85 edition.

\bibitem[\protect\citeauthoryear{}{{Bureau of Public Roads}}{1964}]{BPR}
{Bureau of Public Roads} (1964).
\newblock {\em Traffic Assignment Manual for Application with a Large, High
  Speed Computer}, Vol.~37.
\newblock US Department of Commerce, Bureau of Public Roads, Office of
  Planning.

\bibitem[\protect\citeauthoryear{}{{Bureau of Transportation
  Statistics}}{2020}]{faf4}
{Bureau of Transportation Statistics} (2020).
\newblock ``Freight analysis framework version 4.''\ {\em
  \url{https://faf.ornl.gov/fafweb/}}\ Accessed: 2020-07-25.

\bibitem[\protect\citeauthoryear{}{{Cambridge
  Systematics}}{2012}]{Cambridgesystematics2012}
{Cambridge Systematics} (2012).
\newblock ``Task 8.3: Analysis of freight rail electrification in the scag
  region.''\ {\em
  \url{http://www.freightworks.org/DocumentLibrary/CRGMSAIS\%20-\%20Analysis\%20of\%20Freight\%20Rail\%20Electrification\%20in\%20the\%20SCAG\%20Region.pdf}}\
  Accessed: 2020-07-27.

\bibitem[\protect\citeauthoryear{}{{Center for Transportation
  Analysis}}{2014}]{CTA}
{Center for Transportation Analysis} (2014).
\newblock ``{CTA Railroad Network}\ Accessed: 2018-10-03.

\bibitem[\protect\citeauthoryear{}{Chen and Alfa}{1991}]{Chen}
Chen, M. and Alfa, A.~S. (1991).
\newblock ``A network design algorithm using a stochastic incremental traffic
  assignment approach.''\ {\em Transportation Science}, 25(3), 215--224.

\bibitem[\protect\citeauthoryear{}{Clarke}{1995}]{clarke1995}
Clarke, D.~B. (1995).
\newblock ``An examination of railroad capacity and its implications for
  rail-highway intermodal transportation.''\ Ph.D. thesis, University of
  Tennessee - Knoxville,

\bibitem[\protect\citeauthoryear{}{Dafermos}{1971}]{dafermos1971extended}
Dafermos, S.~C. (1971).
\newblock ``An extended traffic assignment model with applications to two-way
  traffic.''\ {\em Transportation Science}, 5(4), 366--389.

\bibitem[\protect\citeauthoryear{}{Dafermos}{1972}]{dafermos1972traffic}
Dafermos, S.~C. (1972).
\newblock ``The traffic assignment problem for multiclass-user transportation
  networks.''\ {\em Transportation science}, 6(1), 73--87.

\bibitem[\protect\citeauthoryear{}{Dial}{2006}]{Dial}
Dial, R.~B. (2006).
\newblock ``A path-based user-equilibrium traffic assignment algorithm that
  obviates path storage and enumeration.''\ {\em Transportation Research Part
  B: Methodological}, 40(10), 917--936.

\bibitem[\protect\citeauthoryear{}{Ditmeyer et~al.\@}{1985}]{Ditmeyer1985}
Ditmeyer, S., Martin, J., Olson, P., Rister, M., Ross, B., Schmidt, J., and
  Sjokvist, E. (1985).
\newblock ``Railway electrification and railway productivity: A study
  report.''\ {\em Transportation Research Record}, 1029.

\bibitem[\protect\citeauthoryear{}{Drezner and Wesolowsky}{1997}]{Drezner}
Drezner, Z. and Wesolowsky, G.~O. (1997).
\newblock ``Selecting an optimum configuration of one-way and two-way
  routes.''\ {\em Transportation Science}, 31(4), 386--394.

\bibitem[\protect\citeauthoryear{}{Farahani et~al.\@}{2013}]{Farahani2014}
Farahani, R.~Z., Miandoabchi, E., Szeto, W.~Y., and Rashidi, H. (2013).
\newblock ``A review of urban transportation network design problems.''\ {\em
  European Journal of Operational Research}, 229(2).

\bibitem[\protect\citeauthoryear{}{{Fritz, S.G.}}{2000}]{fritz}
{Fritz, S.G.} (2000).
\newblock ``{Diesel fuel effects on locomotive exhaust emissions}.''\ {\em
  Report No. 08.02062}, {Southwest Research Institute, San Antonio, TX}.

\bibitem[\protect\citeauthoryear{}{Gao et~al.\@}{2005}]{Gao2005}
Gao, Z., Wu, J., and Sun, H. (2005).
\newblock ``Solution algorithm for the bi-level discrete network design
  problem.''\ {\em Transportation Research Part B: Methodological}, 39(6),
  479--495.

\bibitem[\protect\citeauthoryear{}{Gartner}{1980}]{gartner1980optimal}
Gartner, N.~H. (1980).
\newblock ``Optimal traffic assignment with elastic demands: a review part ii.
  algorithmic approaches.''\ {\em Transportation Science}, 14(2), 192--208.

\bibitem[\protect\citeauthoryear{}{Gattuso and
  Restuccia}{2014}]{gattuso2014tool}
Gattuso, D. and Restuccia, A. (2014).
\newblock ``A tool for railway transport cost evaluation.''\ {\em
  Procedia-Social and Behavioral Sciences}, 111, 549--558.

\bibitem[\protect\citeauthoryear{}{Gillespie and Hayes}{2003}]{AREMA2003}
Gillespie, A.~J. and Hayes, H.~I. (2003).
\newblock ``Practical guide to railway engineering.''\ {\em AREMA committee-24,
  USA}.

\bibitem[\protect\citeauthoryear{}{Hay}{1982}]{hay1982railroad}
Hay, W.~W. (1982).
\newblock {\em Railroad engineering}, Vol.~1.
\newblock John Wiley \& Sons.

\bibitem[\protect\citeauthoryear{}{Iliopoulou et~al.\@}{2019}]{Iliopoulou2019}
Iliopoulou, C., Kepaptsoglou, K., and Vlahogianni, E. (2019).
\newblock ``Metaheuristics for the transit route network design problem: a
  review and comparative analysis.''\ {\em Public Transport}, 11(3), 487--521.

\bibitem[\protect\citeauthoryear{}{Jayakrishnan et~al.\@}{1994}]{jaya}
Jayakrishnan, R., Tsai, W.~T., Prashker, J.~N., and Rajadhyaksha, S. (1994).
\newblock ``A faster path-based algorithm for traffic assignment.

\bibitem[\protect\citeauthoryear{}{Katoch et~al.\@}{2021}]{katoch2021review}
Katoch, S., Chauhan, S.~S., and Kumar, V. (2021).
\newblock ``A review on genetic algorithm: past, present, and future.''\ {\em
  Multimedia Tools and Applications}, 80(5), 8091--8126.

\bibitem[\protect\citeauthoryear{}{Kneschke}{1986}]{Kneschke1986}
Kneschke, T. (1986).
\newblock ``Simple method for determination of substation spacing for ac and dc
  electrification systems.''\ {\em IEEE Transactions on Industry Applications},
  IA-22(4).

\bibitem[\protect\citeauthoryear{}{Lawrence et~al.\@}{2019}]{chineserail}
Lawrence, M., Bullock, R., and Ziming, L. (2019).
\newblock ``{China's High-Speed Rail Development}\ Accessed: 2021-10-20.

\bibitem[\protect\citeauthoryear{}{Leblanc}{1975}]{LeBlanc}
Leblanc, L.~J. (1975).
\newblock ``An algorithm for the discrete network design problem.''\ {\em
  Transportation Science}, 9(3), 183--199.

\bibitem[\protect\citeauthoryear{}{Lo and Szeto}{2004}]{Lo2004}
Lo, H.~K. and Szeto, W. (2004).
\newblock ``Planning transport network improvements over time.''\ {\em Urban
  and regional transportation modeling: Essays in honor of David Boyce},
  157--176.

\bibitem[\protect\citeauthoryear{}{Lo and Szeto}{2009}]{Lo2009}
Lo, H.~K. and Szeto, W. (2009).
\newblock ``Time-dependent transport network design under cost-recovery.''\
  {\em Transportation Research Part B: Methodological}, 43(1), 142--158.

\bibitem[\protect\citeauthoryear{}{Long et~al.\@}{2010}]{Long}
Long, J., Gao, Z., Zhang, H., and Szeto, W.~Y. (2010).
\newblock ``A turning restriction design problem in urban road networks.''\
  {\em European Journal of Operational Research}, 206(3), 569--578.

\bibitem[\protect\citeauthoryear{}{Meng et~al.\@}{2001}]{Meng2001}
Meng, Q., Yang, H., and Bell, M.~G. (2001).
\newblock ``An equivalent continuously differentiable model and a locally
  convergent algorithm for the continuous network design problem.''\ {\em
  Transportation Research Part B: Methodological}, 35(1), 83--105.

\bibitem[\protect\citeauthoryear{}{{Ministry of Railways
  (India)}}{2020}]{indianrail}
{Ministry of Railways (India)} (2020).
\newblock ``{Railway Electrification}\ Accessed: 2021-10-20.

\bibitem[\protect\citeauthoryear{}{Mishra et~al.\@}{2016}]{Mishra2016}
Mishra, S., Kumar, A., Golias, M.~M., Welch, T., Taghizad, H., and Haque, K.
  (2016).
\newblock ``Transportation investment decision making for medium to large
  transportation networks.''\ {\em Transportation in Developing Economies}.

\bibitem[\protect\citeauthoryear{}{Nektalova}{2008}]{Nektalova}
Nektalova, T. (2008).
\newblock ``\textit{Energy Density of Diesel Fuel}\ Accessed: 2020-07-30.

\bibitem[\protect\citeauthoryear{}{Patil}{2020}]{github2}
Patil, P. (2020).
\newblock ``Rail electrification data repository.''\ {\em
  \url{https://github.com/PriyadarshanPatil/Rail_electrification_GA}}\
  Accessed: 2020-07-25.

\bibitem[\protect\citeauthoryear{}{Patil and Boyles}{2022}]{stap}
Patil, P.~N. and Boyles, S.~D. (2022).
\newblock ``A fresh look at symmetric traffic assignment and algorithm
  convergence.''\ {\em Proc., 101th Annual Meeting}, Transportation Research
  Board.
\newblock Accepted.

\bibitem[\protect\citeauthoryear{}{Patriksson}{1994}]{Patriksson}
Patriksson, M. (1994).
\newblock {\em The traffic assignment problem: models and methods}.
\newblock VSP.

\bibitem[\protect\citeauthoryear{}{Perederieieva
  et~al.\@}{2015}]{perederieieva2015framework}
Perederieieva, O., Ehrgott, M., Raith, A., and Wang, J.~Y. (2015).
\newblock ``A framework for and empirical study of algorithms for traffic
  assignment.''\ {\em Computers \& Operations Research}, 54, 90--107.

\bibitem[\protect\citeauthoryear{}{Prager}{1954}]{prager1954problems}
Prager, W. (1954).
\newblock {\em Problems of traffic and transportation}.

\bibitem[\protect\citeauthoryear{}{{RailTEC}}{2016}]{RailTEC2016}
{RailTEC} (2016).
\newblock ``Transitioning to a zero or near-zero emission line-haul freight
  rail system in california: Operational and economic considerations.''\ {\em
  \url{https://ww3.arb.ca.gov/railyard/docs/uoi_rpt_06222016.pdf}}\ Accessed:
  2020-07-27.

\bibitem[\protect\citeauthoryear{}{Schwarm}{1977}]{schwarm1977}
Schwarm, E.~G. (1977).
\newblock ``Capital and maintenance costs for fixed railroad electrification
  facilities.''\ {\em Transportation Research Board Special Report}, 180.

\bibitem[\protect\citeauthoryear{}{Szeto et~al.\@}{2010}]{Szeto2010}
Szeto, W.~Y., Jaber, X., and O’Mahony, M. (2010).
\newblock ``Time-dependent discrete network design frameworks considering land
  use.''\ {\em Computer-Aided Civil and Infrastructure Engineering}, 25(6).

\bibitem[\protect\citeauthoryear{}{Szeto and Lo}{2006}]{Szeto2006}
Szeto, W.~Y. and Lo, H.~K. (2006).
\newblock ``Transportation network improvement and tolling strategies: the
  issue of intergeneration equity.''\ {\em Transportation Research Part A:
  Policy and Practice}, 40(3), 227--243.

\bibitem[\protect\citeauthoryear{}{Szeto and Lo}{2008}]{Szeto2008}
Szeto, W.~Y. and Lo, H.~K. (2008).
\newblock ``Time-dependent transport network improvement and tolling
  strategies.''\ {\em Transportation Research Part A: Policy and Practice},
  42(2), 376--391.

\bibitem[\protect\citeauthoryear{}{Uddin and Huynh}{2015}]{uddin2015}
Uddin, M.~M. and Huynh, N. (2015).
\newblock ``Freight traffic assignment methodology for large-scale road--rail
  intermodal networks.''\ {\em Transportation Research Record}, 2477(1),
  50--57.

\bibitem[\protect\citeauthoryear{}{{United States Department of Transportation
  - Federal Railroad Administration}}{2015}]{USDOT}
{United States Department of Transportation - Federal Railroad Administration}
  (2015).
\newblock ``Cost-benefit analysis of rail electrification for next generation
  freight and passenger rail transportation.''\ {\em
  \url{https://cms8.fra.dot.gov/sites/fra.dot.gov/files/fra_net/19061/Cost\%20Benefit\%20Analysis\%20of\%20Rail\%20Electrification.pdf}}.

\bibitem[\protect\citeauthoryear{}{{U.S. Bureau of Labor
  Statistics}}{2020}]{bls}
{U.S. Bureau of Labor Statistics} (2020).
\newblock ``Producer price indexes.''\ {\em \url{https://www.bls.gov/ppi}}.

\bibitem[\protect\citeauthoryear{}{{U.S. Department of Transportation - Federal
  Railroad Administration}}{2019}]{USDOT2}
{U.S. Department of Transportation - Federal Railroad Administration} (2019).
\newblock ``Federal-state partnership for state of good repair grant program
  (fy 2019).''\ {\em
  \url{https://railroads.dot.gov/grants-loans/competitive-discretionary-grant-programs/federal-state-partnership-state-good-repair-0}}\
  Accessed: 2020-05-25.

\bibitem[\protect\citeauthoryear{}{{U.S. Geological Survey}}{2007}]{USGS}
{U.S. Geological Survey} (2007).
\newblock ``North america elevation 1-kilometer resolution grid.''\ {\em
  \url{https://www.sciencebase.gov/catalog/item/4fb5495ee4b04cb937751d6d}}.

\bibitem[\protect\citeauthoryear{}{Van~Vuren and
  Watling}{1991}]{van1991multiple}
Van~Vuren, T. and Watling, D. (1991).
\newblock ``A multiple user class assignment model for route guidance.''\ {\em
  Transportation research record},  22--22.

\bibitem[\protect\citeauthoryear{}{Walthall}{2019}]{walthall}
Walthall, R. (2019).
\newblock ``{Rail electrification’s potential for emissions abatement in the
  freight industry : a case study of a transcontinental rail corridor}.''\ M.S.
  thesis, The University of Texas at Austin, ,
  $<$https://hdl.handle.net/2152/80502$>$.

\bibitem[\protect\citeauthoryear{}{Wang et~al.\@}{2018}]{wang2018modeling}
Wang, H., Nozick, L., Xu, N., and Gearhart, J. (2018).
\newblock ``Modeling ocean, rail, and truck transportation flows to support
  policy analysis.''\ {\em Maritime Economics \& Logistics}, 20(3), 327--357.

\bibitem[\protect\citeauthoryear{}{Wardrop and Whitehead}{1952}]{Wardrop}
Wardrop, J.~G. and Whitehead, J.~I. (1952).
\newblock ``Correspondence. some theoretical aspects of road traffic
  research..''\ {\em Proceedings of the Institution of Civil Engineers}, 1(5),
  767--768.

\bibitem[\protect\citeauthoryear{}{Whitford}{1981}]{Whitford1981}
Whitford, R.~K. (1981).
\newblock {\em Railroad Electrification: An Alternative for Petroleum Saving}.
\newblock Automotive Transportation Center, Purdue University.

\end{thebibliography}
%

%FOR PD: For template/PDF submission instructions, see the Template pdf in this repo

\end{document}